\documentclass[psamsfonts,reqno]{amsart}
\usepackage{amssymb,eucal}

\numberwithin{equation}{section}

\newcommand{\pup}[1]{\textup{(}#1\textup{)}}

\theoremstyle{plain}
\newtheorem{lemma}{Lemma}[section]
\newtheorem{theorem}[lemma]{Theorem}
\newtheorem{proposition}[lemma]{Proposition}
\newtheorem{corollary}[lemma]{Corollary}

\newtheorem*{stat}{\name}
\newcommand{\name}{testing}

\theoremstyle{definition}
\newtheorem{definition}[lemma]{Definition}

\newtheorem{problem}{Problem}

\theoremstyle{remark}
\newtheorem{remark}[lemma]{Remark}
\newtheorem{notation}[lemma]{Notation}

\newenvironment{all}[1]{\renewcommand{\name}{#1}\begin{stat}}
                        {\end{stat}}

\newcommand{\qedc}{{\qed}~{\rm Claim~{\theclaim}.}}

\DeclareMathOperator{\sr}{sr}
\DeclareMathOperator{\SR}{SR}
\newcommand{\URPSR}{\mathrm{URP_{sr}}}
\newcommand{\URPSRp}{\mathrm{URP}^{+}_{\mathrm{sr}}}
\newcommand{\jirr}{join-ir\-re\-duc\-i\-ble}

\newcommand{\jsd}{join-sem\-i\-dis\-trib\-u\-tive}
\newcommand{\jsdy}{join-sem\-i\-dis\-trib\-u\-tiv\-i\-ty}

\newcommand{\set}[1]{\left\{#1\right\}}
\newcommand{\setm}[2]{\set{{#1}\mid{#2}}}
\newcommand{\seq}[1]{\left\langle{#1}\right\rangle}
\newcommand{\seqm}[2]{\seq{{#1}\mid{#2}}}

\DeclareMathOperator{\J}{J}

\DeclareMathOperator{\Conc}{Con_c}
\DeclareMathOperator{\Con}{Con}
\DeclareMathOperator{\Dim}{Dim}

\newcommand{\res}{\mathbin{\restriction}}
\newcommand{\es}{\varnothing}

\newcommand{\three}{\boldsymbol{3}}

\newcommand{\diat}[1]{{#1}^{[\wedge 2]}}
\newcommand{\jz}{\ensuremath{(\vee,0)}}
\newcommand{\jzu}{\ensuremath{(\vee,0,1)}}
\newcommand{\jh}{join-ho\-mo\-mor\-phism}
\newcommand{\jzh}{\jz-ho\-mo\-mor\-phism}

\newcommand{\jze}{\jz-em\-bed\-ding}

\newcommand{\jzs}{\jz-se\-mi\-lattice}
\newcommand{\js}{join-se\-mi\-lattice}
\newcommand{\jzus}{\jzu-se\-mi\-lattice}

\newcommand{\cm}{commutative monoid}
\newcommand{\poag}{partially ordered abelian group}

\newcommand{\NN}{\mathbb{N}}
\newcommand{\ZZ}{\mathbb{Z}}

\newcommand{\BB}{\mathcal{B}}
\newcommand{\II}{\mathcal{I}}
\newcommand{\FF}{\mathcal{F}}
\newcommand{\fin}{\mathsf{fin}}
\newcommand{\dnw}{\mathbin{\downarrow}}
\newcommand{\upw}{\mathbin{\uparrow}}

\newcommand{\cC}{\boldsymbol{\mathcal{C}}}

\newcommand{\la}{\boldsymbol{a}}
\newcommand{\lb}{\boldsymbol{b}}

\renewcommand{\le}{\boldsymbol{e}}

\newcommand{\lA}{\boldsymbol{A}}
\newcommand{\lB}{\boldsymbol{B}}

\begin{document}

\title[Semilattices of finitely generated ideals]%
{Semilattices of finitely generated ideals of exchange rings
with finite stable rank}

\author{F. Wehrung}
 \address{CNRS, UMR 6139\\
          Universit\'e de Caen, Campus II\\
          D\'epartement de Math\'ematiques\\
          B.P. 5186\\
          14032 CAEN Cedex\\
          FRANCE}
 \email{wehrung@math.unicaen.fr}
 \urladdr{http://www.math.unicaen.fr/\~{}wehrung}
\date{\today}

\begin{abstract}
We find a distributive \jzus\ $S_{\omega_1}$ of size $\aleph_1$ that
is not isomorphic to the maximal semilattice quotient of any Riesz monoid
endowed with an order-unit of finite stable rank. We thus obtain solutions
to various open problems in ring theory and in lattice theory. In
particular:
\begin{itemize}
\item[---] There is no exchange ring (thus, no von~Neumann regular ring
and no C*-algebra of real rank zero) with finite stable rank whose
semilattice of finitely generated, idempotent-generated two-sided ideals is
isomorphic to~$S_{\omega_1}$.

\item[---] There is no locally finite, modular lattice whose semilattice of
finitely generated congruences is isomorphic to $S_{\omega_1}$.
\end{itemize}
These results are established by constructing an infinitary statement,
denoted here by $\URPSR$, that holds in the maximal semilattice quotient of
every Riesz monoid endowed with an order-unit of finite stable rank,
but not in the semilattice~$S_{\omega_1}$.
\end{abstract}

\maketitle

\section{Introduction}\label{S:Intro}

The present work originates in the representation problem of distributive
\jzs s as semilattices of \emph{compact} (i.e., finitely generated) ideals
of von~Neumann regular rings, see Section~13 in the chapter ``Recent
Developments'' in  \cite{Gvnrr}, or the survey paper \cite{GoWe1}.
It is known that the original form of this problem has a negative solution,
obtained by the author in \cite{Wehr98a,Wehr99}: \emph{There exists a
distributive \jzs\ that is not isomorphic to the semilattice of all compact
ideals of any von Neumann regular ring}. However, the example thus obtained
has size $\aleph_2$. This cardinality bound turns out to be \emph{optimal}.
More specifically: 
\begin{itemize}
\item[---] The main result in \cite{Wehr00} states that every
distributive \jzs\ of size at most $\aleph_1$ is isomorphic to the
semilattice of all compact ideals of some von Neumann regular (but
not locally matricial) algebra.

\item[---] By a well-known unpublished note by G.\,M. Bergman
\cite{Berg86}, every \emph{countable} distributive \jzs\ is
isomorphic to the semilattice of all compact ideals of some locally
matricial (thus unit-regular) algebra. Two other proofs of this result
are presented in \cite{GoWe1}, one of them using the following equivalent
form: \emph{Every countable distributive \jzs\ is isomorphic to the maximal
semilattice quotient $\nabla(G^+)$ of the positive cone $G^+$ of some
dimension group $G$}.

\end{itemize}

Although the cardinality gap seems, at first glance, to be filled by the
above works, the stronger representability result obtained by Bergman
(locally matricial algebras are very special sorts of regular rings) leads
to the question whether similar results hold for distributive semilattices
of cardinality $\aleph_1$.

It turns out that this problem has already been considered by a number of
researchers.

\begin{all}{The $\aleph_1$ nabla problem}
Is every distributive \jzs\ of size $\aleph_1$ isomorphic to the
semilattice of all compact ideals of some locally matricial ring?
\end{all}

To our knowledge, this problem was first published as Problem~10.1
in \cite{GoWe1}. It is also equivalent to Problem~3 in the list of
twelve open problems concluding the survey paper \cite{TuWe2}.

A very important step towards a solution of the $\aleph_1$ nabla problem is
obtained by P. R\r{u}\v{z}i\v{c}ka in \cite{Ruzi1}, where it is proved,
\emph{via} an ingenious construction, that
\emph{there exists a distributive \jzs\ that is not isomorphic to the
maximal semilattice quotient of the positive cone of any dimension group}.
In size $\aleph_1$, the bridge between locally matricial rings and dimension
groups is obtained by results in~\cite{GoHa86}. However,
P.~R\r{u}\v{z}i\v{c}ka's counterexample has size~$\aleph_2$, thus it
does not imply \emph{a priori} a negative solution to the $\aleph_1$ nabla
problem.

In the present paper, we obtain a full negative solution to the $\aleph_1$
nabla problem, see Theorems~\ref{T:URPH} and \ref{T:SnoURPH}. Our
counterexample $S_{\omega_1}$, introduced at the beginning of
Section~\ref{S:Exple}, is very easy to
describe. The idea underlying its proof is a combination of the ideas of
\cite[Example~11.1]{TuWe1} and \cite[Section~3]{Ruzi1}, and it is based on
the construction of a certain infinitary semilattice-theoretical sentence
$\URPSR$, a so-called \emph{uniform
refinement property}, that we present in Section~\ref{S:URPH}. This
property $\URPSR$ holds for the maximal semilattice quotient $\nabla(M)$ of
any Riesz monoid $M$ endowed with an order-unit of finite stable rank
(Theorem~\ref{T:URPH}). On the other hand, it is a near triviality to
verify that $S_{\omega_1}$ does not satisfy $\URPSR$
(Theorem~\ref{T:SnoURPH}); hence
$S_{\omega_1}$ is a counterexample to the $\aleph_1$ nabla problem. The
negative property established for this example is stronger than the one
considered in
\cite{Ruzi1}, both from set-theoretical (cardinality) and algebraic
(strength of the assumptions on $M$) viewpoints.

In Section~\ref{S:Coroll}, we describe some consequences of
Theorems~\ref{T:URPH} and \ref{T:SnoURPH}. The semilattice counterexample
$S_{\omega_1}$ of Section~\ref{S:Exple} is not
isomorphic to the semilattice of all compact ideals of any von~Neumann
regular ring with finite stable rank (Corollary~\ref{C:CXring}), in
particular, $S_{\omega_1}$ is not isomorphic to the semilattice of all
compact ideals of any unit-regular ring. An extension of this result to
\emph{exchange rings} (a more general class including von~Neumann regular
rings and C*-algebras of real rank zero) is proposed in
Corollary~\ref{C:CXExring}. Furthermore, $S_{\omega_1}$ is not isomorphic
to the semilattice of all compact congruences of any modular lattice which
is locally finite, or, more generally, of locally finite length
(Corollary~\ref{C:LocFinLatt}). We also prove that $S_{\omega_1}$ is not
isomorphic to the semilattice of all compact congruences of any locally
finite ``lower bounded'' lattice (Proposition~\ref{P:SnoCLB}), however, we
also prove an apparently much stronger statement in
Corollary~\ref{C:No3Con}.

In Section~\ref{S:PresURP}, we prove a few results, either positive or
negative, about the preservation of $\URPSR$ and related properties under
direct limits of countable sequences.

We conclude the paper in Section~\ref{S:Pbs}, with a few open problems.

\section{Basic concepts}\label{S:Basic}

We mainly follow the notation and terminology in K.\,R.
Goodearl's monograph~\cite{Gpoag}. For elements $a$, $b_1$, \dots, $b_n$ in
some partially ordered set $P$, let $a\leq b_1,\dots,b_n$ abbreviate the
conjunction of the statements $a\leq b_i$, for $1\leq i\leq n$. We put
 \[
 \dnw X=\setm{p\in P}{\exists x\in X\text{ such that }p\leq x},
 \text{ for all }X\subseteq P,
 \]
and we say that $X$ is \emph{cofinal} in $P$, if $\dnw X=P$. We also put
$\dnw p=\dnw\set{p}$, for all $p\in P$. Let $\upw X$, being `coinitial', and
$\upw p$ be defined dually.

For a \poag\ $G$, we denote by $G^+$ the
positive cone of~$G$. We put $\NN=\ZZ^+\setminus\set{0}$. We denote
by $\omega$ (resp., $\omega_1$) the first infinite (resp., uncountable)
ordinal. Following the usual set-theoretical convention, every
ordinal~$\alpha$ is the set of all ordinals less than $\alpha$. For
example, $\omega\setminus n=\setm{k<\omega}{n\leq k}$, for all $n<\omega$.

All our monoids will be denoted additively. Every \cm\ $(M,+,0)$ will be
endowed with its \emph{algebraic} quasi-ordering, defined by $x\leq y$
if{f} there exists $z\in M$ such that $x+z=y$.
We define binary relations $\propto$ and $\asymp$ on $M$ by
 \begin{align*}
 a\propto b&\Longleftrightarrow\exists n\in\NN\text{ such that }a\leq nb,\\
 a\asymp b&\Longleftrightarrow a\propto b\text{ and }b\propto a,
 \end{align*}
for all $a$, $b\in M$. Thus $\asymp$ is a monoid congruence of $M$, the
quotient $\nabla(M)=M/{\asymp}$ is a semilattice, usually called the
\emph{maximal semilattice quotient of} $M$. We denote by $[a]$ the
$\asymp$-equivalence class of $a$, for all $a\in M$, and we call the map
$M\twoheadrightarrow\nabla(M)$, $a\mapsto[a]$ the \emph{canonical
projection} from $M$ onto $\nabla(M)$; see also \cite{GoWe1,Ruzi1}.

An \emph{ideal} of $M$ is a subset $I$ of $M$ such that $0\in I$ and $x+y\in
I$ if{f} $x$, $y\in I$, for all $x$, $y\in M$. An element $e\in M$ is an
\emph{order-unit} of $M$, if $x\propto e$ holds for all $x\in M$.

We say that $M$ is \emph{cancellative}, if $a+c=b+c$ implies that $a=b$,
for all $a$, $b$, $c\in M$. We say that $M$ is \emph{strongly separative},
if $a+b=2b$ implies that $a=b$, for all $a$, $b\in M$. Of course, every
cancellative \cm\ is strongly separative.

We say that $M$ is a \emph{refinement monoid}, if $a_0+a_1=b_0+b_1$ in $M$
implies the existence of $c_{i,j}\in M$, for $i$, $j<2$, such that
$a_i=c_{i,0}+c_{i,1}$ and $b_i=c_{0,i}+c_{1,i}$, for all $i<2$. We say that
$M$ is a \emph{Riesz monoid}, if $c\leq a+b$ in $M$ implies that there are
$a'\leq a$ and $b'\leq b$ in $M$ such that $c=a'+b'$. Every refinement
monoid is a Riesz monoid. The converse does not hold; however, the two
definitions are equivalent for~$M$ a semilattice. We call a semilattice
satisfying these equivalent conditions \emph{distributive}, see
\cite{Grat98}.

\begin{remark}\label{Rk:AntQuot}
Let $M$ be a \cm, define an equivalence relation $\equiv$ on~$M$ by
$x\equiv y$ if $x\leq y\leq x$, for all $x$, $y\in M$. Then $\equiv$ is a
monoid congruence, and we call the quotient $M/{\equiv}$ the \emph{maximal
antisymmetric quotient} of $M$. A sophisticated counterexample by C.
Moreira dos Santos (see \cite{CMor2}) shows that \emph{even for a strongly
separative refinement monoid $M$, the maximal antisymmetric quotient
$M/{\equiv}$ may not have refinement}. On the other hand, it is obvious
that if $M$ is a Riesz monoid, then so is $M/{\equiv}$. Observe that
$\nabla(M)\cong\nabla(M/{\equiv})$.
\end{remark}

We say that a \poag\ $G$ is an \emph{interpolation group}
(see~\cite{Gpoag}), if its positive cone $G^+$ is a refinement monoid. We
say that $G$ is \emph{directed}, if $G=G^++(-G^+)$, and
\emph{unperforated}, if  $mx\geq0$ implies that $x\geq0$, for all $m\in\NN$
and all $x\in G$. A \emph{dimension group} is an unperforated, directed
interpolation group, see also \cite{EHS80}.

In order to conveniently formulate our next lemma, we shall introduce some
notation.

\begin{notation}\label{Not:2Inters}
For a \cm\ $M$ and a family $(a_i)_{i\in I}$ of elements of~$M$,
we denote by $\diat{\seqm{a_i}{i\in I}}$ the set of all elements of $M$ of
the form
 \[
 \sum\bigl(x_p\mid p\subseteq I',\ |p|=2),
 \]
for some finite subset $I'$ of $I$ and
with $x_{\set{i,j}}\leq a_i,a_j$, for all distinct $i$, $j\in I'$.
We also put $\diat{A}=\diat{\seqm{a}{a\in A}}$, for any $A\subseteq M$.
\end{notation}

\begin{lemma}\label{L:AddOrth}
Let $M$ be a Riesz monoid, let $n\in\NN$, and let $a_0$, \dots,
$a_{n-1}$, $b\in M$. If $a_i\leq b$, for all $i<n$, then there exists
$x\in\diat{\seqm{a_i}{i<n}}$ such that $\sum_{i<n}a_i\leq b+x$.
\end{lemma}

\begin{proof}
By induction on $n$. The case $n=1$ is trivial. For $n=2$, let $a_0$,
$a_1\leq b$ in~$M$. There exists $t\in M$ such that $a_1+t=b$. By
$a_0\leq a_1+t$ and since~$M$ is a Riesz monoid, there exists
$x\leq a_0,a_1$ such that $a_0\leq x+t$. Therefore,
$a_0+a_1\leq a_1+t+x=b+x$.

Suppose the result established for some $n$. Let $a_0$, $a_1$, \dots, $a_n$,
$b\in M$ such that $a_i\leq b$, for all $i\leq n$. Put $a=\sum_{i\leq n}a_i$
and $a'=\sum_{i<n}a_i$. Also put $A_k=\diat{\seqm{a_i}{i<k}}$, for all
$k\leq n+1$. By the induction hypothesis, there exists $x\in A_n$ such that
$a'\leq b+x$. Since $a_n\leq b+x$ and $M$ is a Riesz monoid, there exists
$y\leq a',a_n$ such that $a\leq b+x+y$. {}From $y\leq a'=\sum_{i<n}a_i$
and the assumption that $M$ is a Riesz monoid, it follows that there are
$a'_0\leq a_0$, \dots, $a'_{n-1}\leq a_{n-1}$ such that $y=\sum_{i<n}a'_i$.
{}From $y\leq a_n$, it follows that $a'_i\leq a_i,a_n$, for all $i<n$. But
$x\in A_n$, thus $x+y\in A_{n+1}$.
\end{proof}

\section{Stable rank in \cm s}\label{S:StRk}

Let $M$ be a \cm, and let $k$ be a positive integer. An element $e\in M$ has
\emph{stable rank at most $k$}, if $ke+a=e+b$ implies that $a\leq b$,
for all $a$, $b\in M$. Of course, this is equivalent to saying that
$ke+a\leq e+b$ implies that $a\leq b$, for all $a$, $b\in M$. This
definition has at its origin a purely ring-theoretical notion, the
\emph{Bass stable rank} of a ring, see P. Ara's survey paper \cite{AraS}. In
particular, an exchange ring $R$ has stable rank at most $k$ if{f} the
isomorphism class of $R$ has stable rank at most~$k$ in the \cm\ $V(R)$ of
isomorphism classes of finitely generated projective right $R$-modules, see
\cite[Theorem~2.2]{AraS}.

For an element $e$ of $M$, we denote by $\sr_M(e)$ the least positive
integer $k$ such that $e$ has stable rank at most $k$ in $M$ if it exists,
$\infty$ otherwise.

In any \cm, $2e+a=e+b$ implies that $2(e+a)=(e+a)+b$, hence, in the presence
of strong separativity, $e+a=b$. In particular, we obtain the
following result.

\begin{proposition}\label{P:StRkStSep}
For any \cm\ $M$, the following statements hold:
\begin{enumerate}
\item If $M$ is cancellative, then every element of $M$ has stable
rank at most $1$.
\item If $M$ is strongly separative, then every element of $M$ has stable
rank at most $2$.
\end{enumerate}
\end{proposition}

We put
 \[
 \SR(M)=\setm{(x,y)\in M\times M}
 {\forall a,b\in M,\ x+a\leq y+b\Rightarrow a\leq b}.
 \]
Obviously, $\SR(M)$ is a submonoid of $M\times M$. Observe that
$\sr_M(e)\leq k$ if{f} $(ke,e)\in\SR(M)$. As an immediate consequence of
these two simple facts, we observe the following.

\begin{lemma}\label{L:sr(e)}
Let $e$ be an element of a \cm\ $M$. Then $\sr_M(ne)\leq\sr_M(e)$, for all
$n\in\NN$.
\end{lemma}

\section{A new uniform refinement property}\label{S:URPH}

The negative solutions to other representation problems, either of dimension
groups or distributive semilattices, introduced in \cite{Wehr98a,Wehr99},
were obtained as counterexamples to special infinitary statements called
\emph{uniform refinement properties}. The reader may consult \cite{TuWe2}
for a discussion about various such statements defined for semilattices.

Our new uniform refinement property is significantly different from all
previously known ones. It is inspired by the proof of the
counterexample in \cite[Example~11.1]{TuWe1}. In order to state it, we find
it convenient to use Notation~\ref{Not:2Inters}.  We say that a partially
ordered set $P$ is \emph{$\aleph_0$-downward directed}, if every at
most countable subset of $P$ lies above some element of $P$. We shall use
the following straightforward lemma several times.

\begin{lemma}\label{L:ExtrCof}
Let $(P_n)_{n<\omega}$ be a sequence of subsets of a $\aleph_0$-downward
directed partially ordered set $P$. If $\bigcup_{n<\omega}P_n=P$, then
there exists $n<\omega$ such that $P_n$ is coinitial in $P$.
\end{lemma}

Now we turn to the most important definition of the paper.

\begin{definition}\label{D:URPS}
Let $S$ be a \js. For an element $e\in S$, we introduce the
following statements:
\begin{itemize}
\item $\URPSR(e)$: for all subsets $A$ and $B$ of $S$ such that
$A$ is uncountable, $B$ is $\aleph_0$-downward directed, and
 \[
 a\leq e\leq a\vee b\text{ for all }(a,b)\in A\times B,
 \]
there exists $a\in\diat A$ such that $e\leq a\vee b$ for all $b\in B$.

\item $\URPSRp(e)$: for all $a_0$, $a_1\in S$
and all $\aleph_0$-downward directed $B\subseteq S$, if
 \[
 e\leq a_i\vee b\text{ for all }i<2\text{ and all }b\in B,
 \]
then there exists $a\leq a_0,a_1$ such that $e\leq a\vee b$ for all
$b\in S$.
\end{itemize}

We say that $S$ satisfies $\URPSR$ (resp., $\URPSRp$), if it satisfies
$\URPSR(e)$ (resp., $\URPSRp(e)$), for every $e\in S$.
\end{definition}

Of course, $\URPSRp$ implies $\URPSR$.

Our main result is the following.

\begin{theorem}\label{T:URPH}
Let $M$ be a Riesz monoid and let $e\in M$. If $e$ has finite stable rank
in~$M$, then the maximal semilattice quotient $\nabla(M)$ satisfies
$\URPSR([e])$.
\end{theorem}

\begin{proof}
Put $S=\nabla(M)$ and $\le=[e]\in S$. Let $\lA$
and $\lB$ be subsets of $S$ such that $\lA$ is uncountable, $\lB$ is
$\aleph_0$-downward directed, and $\la\leq\le\leq\la\vee\lb$ for all
$(\la,\lb)\in\lA\times\lB$. It is convenient to write these sets as
 \[
 \lA=\setm{\la_i}{i\in I}\text{ and }\lB=\setm{\lb_j}{j\in J},
 \]
where $i\mapsto\la_i$ is one-to-one, $J$ is a $\aleph_0$-downward directed
partially ordered set, and $j\mapsto\lb_j$ is order-preserving.
Pick $a_i\in\la_i$ for all $i\in I$ and
$b_j\in\lb_j$ for all $j\in J$. For all $i\in I$, there exists $m_i\in I$
such that $a_i\leq m_ie$. Since $I$ is uncountable, there are $m\in\NN$ and
an uncountable subset $U$ of $I$ such that $m_i=m$ for all $i\in U$. Thus
$a_i\leq me$ for all $i\in U$. But $me$ has finite stable rank in $M$ (see
Lemma~\ref{L:sr(e)}), thus we may replace $e$ by $me$, so that we obtain
 \begin{equation}\label{Eq:exie}
 a_i\leq e,\text{ for all }i\in U.
 \end{equation}
For all $(i,j)\in I\times J$, since $\le\leq\la_i\vee\lb_j$, there exists
$n_{i,j}\in\NN$ such that $e\leq n_{i,j}(a_i+b_j)$. Put
$J_{i,n}=\setm{j\in J}{n_{i,j}=n}$, for all $(i,n)\in I\times\NN$. For all
$i\in I$, the equality $J=\bigcup_{n\in\NN}J_{i,n}$ holds. Thus, by
Lemma~\ref{L:ExtrCof}, there exists $n_i\in\NN$ such that $J_{i,n_i}$ is
coinitial in $J$. Since $U$ is uncountable, there are $n\in\NN$ and an
uncountable subset~$V$ of $U$ such that $n_i=n$ for all $i\in V$. Observe
that, by construction,
 \begin{equation}\label{Eq:eaibj}
 e\leq n(a_i+b_j),\text{ for all }i\in V\text{ and }j\in J_{i,n}.
 \end{equation}
Put $k=\sr_M(e)$.
We pick a subset $X\subset V$ with exactly $nk+1$ elements. By
\eqref{Eq:exie}, the fact that $X$ is a subset of $U$, and
Lemma~\ref{L:AddOrth}, it follows that
 \begin{equation}\label{Eq:aie+a}
 \sum_{i\in X}a_i\leq e+a,\text{ for some }a\in\diat{\seqm{a_i}{i\in X}}.
 \end{equation}
Observe that the element $\la=[a]$ of $S$ belongs to
$\diat{\seqm{\la_i}{i\in X}}$, thus
(since $i\mapsto\la_i$ is one-to-one) to $\diat{\lA}$.
Let $j\in J$. For all $i\in I$, the subset $J_{i,n}=J_{i,n_i}$ is coinitial
in $J$, thus there exists $\varphi(i)\in J_{i,n}$ such that
$\varphi(i)\leq j$. It follows from \eqref{Eq:eaibj} that
$e\leq n(a_i+b_{\varphi(i)})$. Therefore, adding together all those
inequalities for $i\in X$ and using \eqref{Eq:aie+a}, we obtain that
 \[
 (nk+1)e\leq n\sum_{i\in X}a_i+n\sum_{i\in X}b_{\varphi(i)}
 \leq ne+na+n\sum_{i\in X}b_{\varphi(i)},
 \]
thus, since $\sr_M(ne)\leq k$ (see Lemma~\ref{L:sr(e)}),
 \[
 e\leq na+n\sum_{i\in X}b_{\varphi(i)}.
 \]
By applying to this inequality the canonical projection from $M$ onto
$\nabla(M)$ and by using the fact that $\lb_{\varphi(i)}\leq\lb_j$ for all
$i\in X$, we obtain the inequalities
 \[
 \le\leq\la\vee\bigvee_{i\in X}\lb_{\varphi(i)}\leq\la\vee\lb_j.
 \]
Hence the element $\la\in S$ is as required.
\end{proof}

\begin{remark}\label{Rk:URPH}
Although the assumptions made in the statement of Theorem~\ref{T:URPH} are
not the weakest possible, they are probably the weakest that can be stated
conveniently while meeting an application range as wide as possible within
mathematical practice. Nevertheless, since one can never be sure about
future applications, we list here some possible weakenings of the
assumptions of Theorem~\ref{T:URPH} that lead to the same conclusion.
\begin{itemize}
\item[---] Refinement assumption: for all $m\in\NN$ and all $a_0$, \dots,
$a_{m-1}$, $b\in M$, if $a_i\leq b$ for all $i<m$, then there are $n\in\NN$
and $x\in\diat{\seqm{na_i}{i<m}}$ such that $\sum_{i<m}a_i\leq b+x$.

\item[---] Stable rank assumption: for all $m\in\NN$, there exists $k\in\NN$
such that for all $a$, $b\in M$, if $kme+a\leq me+b$, then $a\propto b$.
\end{itemize}
\end{remark}

\begin{corollary}\label{C:URPH}
Let $M$ be a strongly separative Riesz monoid. Then $\nabla(M)$ satisfies
$\URPSR$.
\end{corollary}

\section{The example}\label{S:Exple}

We first recall some constructions used in \cite{Ruzi1}. For an infinite
cardinal $\kappa$, we denote by $\BB_{\kappa}$
the \emph{interval algebra} of $\kappa$, that is, the Boolean subalgebra
of the powerset of $\kappa$ generated by all intervals of $\kappa$.
Then we put
 \begin{align*}
 \II_{\kappa}&=\setm{x\in\BB_{\kappa}}
 {\exists\alpha<\kappa\text{ such that }x\subseteq\alpha},\\
 \FF_{\kappa}&=\setm{x\in\BB_{\kappa}}
 {\exists\alpha<\kappa\text{ such that }\kappa\setminus\alpha\subseteq x},\\
 D_{\kappa}&=\setm{x\subseteq\kappa}
 {\text{either }x\text{ is finite or }x=\kappa}.
 \end{align*}
Observe that $D_{\kappa}$ is a distributive lattice with zero.
We again use the notation $D^-=D\setminus\set{\es}$, and we put
 \[
S_{\kappa}=(\set{\es}\times\II_{\kappa})\cup(D_{\kappa}^-\times\FF_{\kappa})
 \qquad(\text{thus }S_{\kappa}\subseteq D_{\kappa}\times\BB_{\kappa}),
 \]
endowed with its componentwise ordering. This is a particular case of a
construction introduced in \cite{Ruzi1}. Observe that $S_{\kappa}$ is not a
lattice (the meet of two elements of $D_{\kappa}^-$ may be empty), however,
it is a distributive \jzus\ (see \cite[Lemma~3.3]{Ruzi1}). Observe also that
$S_{\kappa}$ has size $\kappa$.

\begin{theorem}\label{T:SnoURPH}
The distributive \jzs\ $S_{\omega_1}$ does not satisfy $\URPSR(e)$, where
we put $e=(\omega_1,\omega_1)$, the largest element of $S_{\omega_1}$.
\end{theorem}

\begin{proof}
We put $A=\setm{a_\xi}{\xi<\omega_1}$ and
$B=\setm{b_\xi}{\xi<\omega_1}$, where we put
 \[
 a_\xi=(\set{\xi},\omega_1)\text{ and }
 b_\xi=(\omega_1,\omega_1\setminus\xi),\text{ for all }\xi<\omega_1.
 \]
It is obvious that $A$ is uncountable, $B$ is
$\aleph_0$-downward directed, and $a\leq e\leq a\vee b$ 
(in fact, $e=a\vee b$) for all $(a,b)\in A\times B$. For any
$a\in\diat{\seqm{a_\xi}{\xi<\omega_1}}$, there exists an ordinal
$\alpha<\omega_1$ such that $a\leq(\es,\alpha)$. Hence
$a\vee b_{\alpha+1}\leq(\omega_1,\omega_1\setminus\set{\alpha})<e$,
which completes the proof.
\end{proof}

Thus, by using Theorem~\ref{T:URPH}, we obtain the following.

\begin{corollary}\label{C:SnoURPH}
Let $M$ be a Riesz monoid in which there is an order-unit of finite
stable rank. Then $S_{\omega_1}$ is not isomorphic to $\nabla(M)$.
\end{corollary}

Also, by using Corollary~\ref{C:URPH}, we obtain the following.

\begin{corollary}\label{C:CXAl1}
There exists no strongly separative refinement monoid $M$ such that
$\nabla(M)\cong S_{\omega_1}$. In particular, there exists no interpolation
group \pup{thus, no dimension group} $G$ such that
$\nabla(G^+)\cong S_{\omega_1}$.
\end{corollary}

\section{Consequences in ring theory and lattice theory}\label{S:Coroll}

In this section we shall reap some consequences of Theorem~\ref{T:URPH}
and Theorem~\ref{T:SnoURPH}, thus answering a few open questions in
ring theory and lattice theory. Throughout this section, we shall denote by
$S_{\omega_1}$ the semilattice constructed in Section~\ref{S:Exple}.

\subsection{Ideal lattices of von Neumann regular rings}\label{Sub:vnrr}
For a von~Neumann regular ring $R$, the monoid $V(R)$ of isomorphism
classes of finitely generated projective right $R$-modules is a refinement
monoid (see \cite[Chapter~2]{Gvnrr}); it is, in addition, \emph{conical},
that is, it satisfies the quasi-identity $x+y=0\Rightarrow x=y=0$.
Furthermore, the lattice of (two-sided) ideals of $R$ is isomorphic to the
lattice of ideals of $V(R)$, see \cite[Proposition~7.3]{GoWe1}.

We say that $R$ is \emph{strongly separative}, if the monoid $V(R)$ is
strongly separative. This notion is the monoid-theoretical translation of
a purely ring-theoretical notion, see \cite{AGPO}. We observe, for
a given von Neumann regular ring $R$, that the following implications hold:
 \[
 \text{locally matricial}\Rightarrow
 \text{unit-regular}\Rightarrow\text{strongly separative}
 \Rightarrow\text{finite stable rank}.
 \]

\begin{corollary}\label{C:CXring}
There is no von Neumann regular ring $R$ with finite stable rank such that
$S_{\omega_1}$ is isomorphic to the semilattice of all compact ideals of
$R$.
\end{corollary}

Corollary~\ref{C:CXring} implies immediately that $S_{\omega_1}$ is not
isomorphic to the semilattice of all compact ideals of any strongly
separative (resp., unit-regular, locally matricial) ring.
By contrast, we recall that every countable distributive \jzs\ is
isomorphic to the semilattice of all compact ideals of some countable
dimensional locally matricial ring (over any given field), see
\cite{Berg86,GoWe1}.

We can even formulate an extension of Corollary~\ref{C:CXring} to more
general rings called \emph{exchange rings}. Recall (see the survey paper
\cite{AraS}) that every von~Neumann regular ring or every C*-algebra with
real rank zero is an exchange ring. For an arbitrary exchange ring $R$, the
canonical map from the ideal lattice of $R$ to the ideal lattice of $V(R)$
is surjective, and the inverse image of a singleton $\set{I}$ is the
interval $[J_0,J_1]$, where $J_0$ is the ideal generated by all idempotents
$e\in R$ such that $[eR]\in I$, and $J_1$ is the intersection of all
primitive ideals containing $J_0$ (see \cite[Teorema~4.1.7]{Pard95}). In
particular, the maximal semilattice quotient of $V(R)$ is isomorphic to the
semilattice of all compact idempotent-generated ideals of $R$. Hence we
obtain the following.

\begin{corollary}\label{C:CXExring}
There is no exchange ring $R$ with finite stable rank such that
$S_{\omega_1}$ is isomorphic to the semilattice of all compact
idempotent-generated ideals of $R$.
\end{corollary}

\subsection{Congruence lattices of modular lattices}\label{Sub:ModLatt}
We say that a lattice $K$ has \emph{locally finite length}, if every
finitely generated sublattice of~$K$ has finite length. In
particular, all locally finite lattices and all direct limits of
finite-dimensional projective geometries have locally finite length.

\begin{corollary}\label{C:LocFinLatt}
There is no modular lattice~$K$ of locally finite length such that the
semilattice $\Conc K$ of all compact congruences of $K$ is isomorphic to
$S_{\omega_1}$.
\end{corollary}

\begin{proof}
Suppose that $S_{\omega_1}\cong\Conc K$, for a modular lattice $K$ of
locally finite length. It follows from \cite[Corollary~2.3]{Wehr98b} that
$S_{\omega_1}$ is isomorphic to $\nabla(M)$, for $M=\Dim K$, the so-called
\emph{dimension monoid} of $K$. However, the dimension monoid of a modular
lattice of finite length is a finitely generated free commutative monoid
(see \cite[Proposition~5.5]{Wehr98b}), and the dimension monoid functor
preserves direct limits (see \cite[Proposition~1.4]{Wehr98b}). Hence, since
$K$ is modular with locally finite length, the monoid $\Dim K$ is a direct
limit of free commutative monoids, and thus it is the positive cone of a
dimension group~$G$. Hence,
 \[
 S_{\omega_1}\cong\nabla(\Dim K)=\nabla(G^+),
 \]
a contradiction by Corollary~\ref{C:CXAl1}.
\end{proof}

In contrast, we recall that every distributive \jzs\ of size at
most~$\aleph_1$ is isomorphic to the semilattice of all compact
congruences of some locally finite relatively complemented lattice with zero
(see \cite{GLWe}), and also isomorphic to the semilattice of all compact
ideals of some von~Neumann regular algebra, and to the semilattice of all
compact congruences of some sectionally complemented, modular (but not
locally finite) lattice (see \cite{Wehr00}).

\subsection{Congruence lattices of lower bounded lattices}\label{Sub:LB}
For lattices $K$ and $L$, a lattice homomorphism $f\colon K\to L$ is
\emph{lower bounded}, if $\setm{x\in K}{a\leq f(x)}$ 
is either empty or has a least element,
for every $a\in L$. As in \cite{AdGo}, we say that a lattice~$L$ is
\emph{lower bounded}, if every homomorphism from a finitely generated
lattice to $L$ is lower bounded. There are many equivalent definitions of
lower boundedness for finite lattices, the simplest of them being that $L$
has as many \jirr\ elements as its congruence lattice $\Con L$, see
\cite[Lemma~2.40]{FJN}. Every lower bounded lattice is
\emph{\jsd}, that is, it satisfies the quasi-identity
 \[
 x\vee y=x\vee z\Longrightarrow x\vee y=x\vee(y\wedge z).
 \]
Both properties of lower boundedness and \jsdy\ are antithetical to
modularity, as, for example, \emph{every \jsd\ modular lattice is
distributive}. This antithesis also explains why results
proved for modular lattices are often worth investigating for \jsd, or lower
bounded, lattices.

Congruence lattices of finite lower bounded lattices are fully
characterized in \cite{GrWe}. In particular,
not every finite distributive lattice is isomorphic to the congruence
lattice of some finite lower bounded lattice, the three-element chain
$\three$ being the simplest such example. On the dimension side, the
following result is a consequence of \cite[Corollary~6.3]{Wehr02}.

\begin{proposition}\label{P:DimLB}
The dimension monoid $\Dim L$ of any finite lower bounded lattice~$L$ is
strongly separative.
\end{proposition}

Since the $\Dim$ functor preserves direct limits, the result of
Proposition~\ref{P:DimLB} extends immediately to \emph{locally finite} lower
bounded lattices. Hence, from the isomorphism $\nabla(\Dim L)\cong\Conc L$
and by Corollary~\ref{C:CXAl1}, we immediately obtain the following.

\begin{proposition}\label{P:SnoCLB}
There is no locally finite, lower bounded lattice $L$ such that
$\Conc L\cong S_{\omega_1}$.
\end{proposition}

On the other hand, a much more striking negative conclusion can be reached
by totally different means, \emph{via} the following result.

\begin{proposition}\label{P:FinCon}
Let $L$ be a lower bounded lattice. If $\Con L$ is finite, then so is~$L$.
\end{proposition}

\begin{proof}
It is proved in \cite{Day79} that every finite lower bounded lattice $L$
satisfies the following \emph{Day-Pudl\'ak-T\r{u}ma property} (DPT):
 \[
 \bigl[\Theta(a',a)=\Theta(b',b)\text{ and }a'<a\text{ and }b'<b\bigr]
 \Rightarrow a\wedge b\nleq a',b',\text{ for all }a,\,a',\,b,\,b'\in L,
 \]
where $\Theta(x,y)$ denotes the congruence of $L$ generated by the pair
$(x,y)$. By \cite[Theorem~2.1]{AdGo}, this result can be extended to
arbitrary lower bounded lattices. Now suppose that $\Con L$ is finite, put
$n=|\Con L|-1$. If there exists a chain $x_0<x_1<\cdots<x_{n+1}$ in $L$,
then there are $i<j$ such that $\Theta(x_i,x_{i+1})=\Theta(x_j,x_{j+1})$.
Thus, by (DPT), $x_{i+1}=x_{i+1}\wedge x_{j+1}\nleq x_j$, a contradiction.
Hence $L$ has length at most~$n$, but it is lower bounded, thus
\jsd, hence, by a classical result of B.~J\'onsson and J.\,E. Kiefer, see
\cite{JoKi} (also \cite[Theorem~5.59]{FJN}), $|L|\leq 2^n$.
\end{proof}

The upper bound $|L|\leq 2^n$ is rather crude. Indeed, since $L$ is finite
lower bounded, there exists (see \cite[Lemma~2.40]{FJN}) a bijection between
the set $\J(L)$ of \jirr\ elements of $L$ and $\J(\Con L)$. In particular,
$|L|\leq 2^{|\J(\Con L)|}$.

\begin{corollary}\label{C:No3Con}
There is no lower bounded lattice $L$ such that $\Con L\cong\three$.
\end{corollary}

\begin{proof}
If $\Con L\cong\three$, then, by Proposition~\ref{P:FinCon}, $L$ is finite.
However, there is no such finite lower bounded lattice $L$, see \cite{GrWe}.
\end{proof}

Observe that no assumption of local finiteness is necessary in
Proposition~\ref{P:FinCon} and Corollary~\ref{C:No3Con}. On the other hand,
unlike $S_{\omega_1}$, the \jzs\ $\three$ is the maximal semilattice
quotient of a strongly separative refinement monoid $M$---in fact, by
Bergman's Theorem, $M$ may be chosen as the positive cone of some dimension
group. See the related Problem~\ref{Pb:IdSSIG}.

\section{Preservation under countable direct limits}
\label{S:PresURP}

The following easy result introduces another class of
distributive semilattices with $\URPSR$, defined in a completely different
way.

\begin{proposition}\label{P:CtbleInSeg}
Let $S$ be a distributive \js\ such that $\dnw a_0\cap\dnw a_1$ has an at
most countable cofinal subset, for all $a_0$, $a_1\in S$. Then $S$
satisfies $\URPSRp$.
\end{proposition}

\begin{proof}
Let $e$, $a_0$, $a_1\in S$ and let $B$ be a $\aleph_0$-downward
directed subset of $S$ such that $e\leq a_i\vee b$ for all $i<2$ and all
$b\in B$. By assumption, there exists an increasing cofinal sequence
$\setm{a^{(n)}}{n<\omega}$ in $\dnw a_0\cap\dnw a_1$. Suppose that for all
$n<\omega$, there exists $b_n\in B$ such that $e\nleq a^{(n)}\vee b_n$. By
assumption, there exists $b\in B$ such that $b\leq b_n$ for all $n<\omega$;
hence $e\nleq a^{(n)}\vee b$ for all $n<\omega$. On the other hand,
$e\leq a_i\vee b$ for all $i<2$, thus, since $S$ is distributive,
$e\leq x\vee b$ for some $x\leq a_0,a_1$. Pick $n<\omega$ such that
$x\leq a^{(n)}$, then $e\leq a^{(n)}\vee b$, a contradiction. So we have
proved that there exists $n<\omega$ such that $e\leq a^{(n)}\vee b$ for all
$b\in B$.
\end{proof}

Observe the following immediate corollary of Proposition~\ref{P:CtbleInSeg}.

\begin{corollary}\label{C:UnDn}
Every direct limit of a countable sequence of distributive lattices and
\jh s satisfies $\URPSRp$ \pup{thus also $\URPSR$}.
\end{corollary}

Compare with \cite[Proposition~2.11]{TuWe2}.

We shall see soon (Proposition~\ref{L:Som1ctblelim}) that the class of all
semilattices satisfying $\URPSR$ is \emph{not} closed under
direct limits of countable sequences. However, the following related
positive result holds.

\begin{proposition}\label{P:CtbleUnURP}
Let $(S_n)_{n<\omega}$ be an increasing
sequence of join-subsemilattices of a \js~$S$ such that
$S=\bigcup_{n<\omega}S_n$. If all $S_n$-s satisfy $\URPSR$ \pup{resp.,
$\URPSRp$}, then so does $S$.
\end{proposition}

\begin{proof}
We provide a proof for $\URPSR$; the proof for $\URPSRp$ is similar.
Let $e\in S$ and let $A$ and $B$ be subsets of $S$, where $A$ is
uncountable, $B$ is $\aleph_0$-downward directed, and $a\leq e\leq a\vee b$
for all $(a,b)\in A\times B$. By Lemma~\ref{L:ExtrCof}, there exists
$n<\omega$ such that $B_n=B\cap S_n$ is coinitial in $B$---in particular,
it is $\aleph_0$-downward directed; furthermore, we may assume that
$A_n=A\cap S_n$ is uncountable. Since $a\leq e\leq a\vee b$ for all
$(a,b)\in A_n\times B_n$ and $S_n$ satisfies $\URPSR$, there exists
$a\in\diat{A_n}$ such that $e\leq a\vee b$ for all $b\in B_n$. So
$a\in\diat{A}$ and, since $B_n$ is coinitial in $B$, the inequality
$e\leq a\vee b$ holds for all $b\in B$.
\end{proof}

Denote by $\cC$ the class of all distributive \jzs s without any decreasing
$\omega_1$-chain. Observe that any nonempty $\aleph_0$-downward directed
subset of any member~$S$ of $\cC$ has a least element. Hence the
following result holds.

\begin{proposition}\label{P:C2URPSRp}
Every member of $\cC$ satisfies $\URPSRp$.
\end{proposition}

Denote by $\cC_\omega$ the class of all direct limits of countable
sequences of members of $\cC$, and by $\fin$ the Fr\'echet filter on
$\omega$. Furthermore, for a set $X$ and a sequence $x=\seqm{x_n}{n<\omega}$
in $X^\omega$, we denote by $[x_n\mid n<\omega]$ the equivalence class of
$x$ modulo $\fin$ in $X^\omega$.

\begin{lemma}\label{L:RPow2DirLim}
Let $S$ and $C$ be distributive \jzs s with $C\in\cC$. If $S$ embeds into
the reduced power $C^\omega/\fin$, then $S$ belongs to $\cC_\omega$.
\end{lemma}

\begin{proof}
Let $f\colon S\hookrightarrow C^\omega/\fin$ be a \jze. Denote by
$S'$ its image. For all $n<\omega$, let
$\pi_n\colon C^{\omega\setminus n}\twoheadrightarrow C^\omega/\fin$
be the map defined by
 \[
 \pi_n(\seqm{x_k}{n\leq k<\omega})=[x_k\mid k<\omega],
 \text{ where we put }x_0=\cdots=x_{n-1}=0.
 \]
It is not hard to verify that $S_n=\pi_n^{-1}[S']$ is a distributive
\jz-subsemilattice of $C^{\omega\setminus n}$. Put
 \[
 \rho_n=f^{-1}\circ\pi_n\res_{S_n}.
 \]
Furthermore, for $m\leq n<\omega$, it is possible to define a \jzh\ 
$\rho_{m,n}\colon S_m\to S_n$ by the rule
$\rho_{m,n}(x)=x\res_{\omega\setminus n}$, for all $x\in S_m$.

It is routine to verify that $S=\varinjlim_{n<\omega}S_n$, with
transition maps $\rho_{m,n}\colon S_m\to S_n$ and limiting maps
$\rho_n\colon S_n\to S$. The semilattice $S_n$ is a distributive
subsemilattice of $C^{\omega\setminus n}$, for all $n<\omega$, and
$C\in\cC$, thus $S_n\in\cC$.
\end{proof}

Now we shall present a construction which proves that the semilattice
$S_{\omega_1}$ of Section~\ref{S:Exple} belongs to $\cC_\omega$. For maps
$f$, $g$ from $\omega$ to $\omega+1$, let $f<^*g$ be the
following statement:
 \[
 f<^*g\Longleftrightarrow\setm{n<\omega}{f(n)\geq g(n)}\text{ is finite}.
 \]
It is well-known that there exists a $\omega_1$-sequence
$(f_\alpha)_{\alpha<\omega_1}$ of maps from $\omega$ to $\omega$ such that
$\alpha<\beta$ implies that $f_\alpha<^*f_\beta$, for all $\alpha$,
$\beta<\omega_1$. Let $f_{\omega_1}\colon\omega\to\omega+1$ denote the
constant function with value $\omega$. Observe that
 \begin{equation}\label{Eq:fa<fb}
 \alpha<\beta\Longrightarrow f_\alpha<^*f_\beta,
 \text{ for all }\alpha,\,\beta\leq\omega_1.
 \end{equation}
Every element $x\in\BB_{\omega_1}$ has a unique \emph{normal form},
 \begin{equation}\label{Eq:NormFormx}
 x=\bigcup_{i<n}[\alpha_i,\beta_i),\text{ where }
 n<\omega\text{ and }\alpha_0<\beta_0<\cdots<\alpha_{n-1}<\beta_{n-1}
 \leq\omega_1.
 \end{equation}
Observe that the $[\alpha_i,\beta_i)$-s are exactly the \emph{maximal
subintervals} of $x$. For $x\in\BB_{\omega_1}$ written in normal form as in
\eqref{Eq:NormFormx}, we put
 \begin{equation}\label{Eq:Defgk}
 g_k(x)=\bigcup_{i<n}[f_{\alpha_i}(k),f_{\beta_i}(k)),\text{ for all }
 k<\omega,
 \end{equation}
where, of course, an interval $[x,y)$ is empty if $x\geq y$.
Observe that for all large enough $k$, the following inequalities hold:
 \[
 f_{\alpha_0}(k)<f_{\beta_0}(k)<\cdots<f_{\alpha_{n-1}}(k)
 <f_{\beta_{n-1}}(k).
 \]
Hence, for such values of $k$, \eqref{Eq:Defgk} is an expression of
$g_k(x)$ in normal form.

For an element $x\in D_{\omega_1}$ and $k<\omega$, we define
$u_k(x)\in D_{\omega}$ by
 \[
 u_k(x)=\begin{cases}
 \setm{f_\alpha(k)}{\alpha\in x},&\text{if }x\text{ is finite},\\
 \omega,&\text{if }x=\omega_1.
 \end{cases}
 \]
Finally, we put
 \[
h(x,y)=[(u_k(x),g_k(y))\mid k<\omega],\text{ for all }(x,y)\in S_{\omega_1}.
 \]
We leave to the reader the easy but somehow tedious proof of the following
lemma. The main reason why it works is that containments between elements
$x$ and $y$ of $S_\kappa$, for an infinite cardinal $\kappa$, can be
expressed by inequalities between the endpoints of the maximal subintervals
of $x$ and $y$. But then, inequalities between ordinals in $\omega_1+1$ can
be ``projected on'' $\omega+1$ by using \eqref{Eq:fa<fb}.

\begin{lemma}\label{L:Embom12om0}
The map $h$ is a \jze\ from $S_{\omega_1}$ into
$(S_{\omega})^\omega/\fin$.
\end{lemma}

Observe that Lemma~\ref{L:Embom12om0} does not hold ``trivially'', in the
following sense. Although it is easy to prove that all the maps $u_k$, for
$k<\omega$, are \jh s, this is not the case for the $g_k$-s.
However, this seemingly irregular behavior disappears ``at the limit'', as
$k$ goes to infinity.

By using Lemmas~\ref{L:RPow2DirLim} and \ref{L:Embom12om0}, we obtain the
following.

\begin{proposition}\label{L:Som1ctblelim}
The semilattice $S_{\omega_1}$ belongs to $\cC_\omega$.
\end{proposition}

By Theorem~\ref{T:SnoURPH}, Proposition~\ref{P:C2URPSRp}, and
Proposition~\ref{L:Som1ctblelim}, the class of all semilattices satisfying
$\URPSR$ is not closed under direct limits of countable sequences. Compare
with Corollary~\ref{C:UnDn} and Proposition~\ref{P:CtbleUnURP}.

\begin{remark}\label{Rk:CtbleLimRepr}
Let $h\colon S_{\omega_1}\hookrightarrow(S_{\omega})^\omega/\fin$ be the
previously constructed \jze. It follows from
Lemma~\ref{L:RPow2DirLim} that $S_{\omega_1}$ is the direct limit of the
lattices $T_n=\pi_n^{-1}h[S_{\omega_1}]$, for $n<\omega$, and all $T_n$-s
belong to $\cC$. Furthermore, the map
$\rho_n=h^{-1}\circ\pi_n\colon T_n\twoheadrightarrow S_{\omega_1}$
is a surjective homomorphism of a very special kind: namely,
 \[
 \rho_0(x)\leq\rho_0(y)\Longleftrightarrow\exists
 u\in S_{\omega}^{(\omega)}
 \text{ such that }x\leq y\vee u,\text{ for all }x,\,y\in T_0,
 \]
where we denote by $S_{\omega}^{(\omega)}$ the ideal of
$(S_{\omega})^{\omega}$ that consists of all sequences with finite support.
Hence,
 \begin{equation}\label{Eq:S1asQuot}
 S_{\omega_1}\cong T_0/S_{\omega}^{(\omega)}.
 \end{equation}
It follows from this that \emph{$T_0$ is not isomorphic to $\nabla(G^+)$,
for any interpolation group~$G$}. Indeed, otherwise, by \eqref{Eq:S1asQuot},
$S_{\omega_1}$ would be isomorphic to $\nabla((G/I)^+)$, where $I$ is
the ideal of $G$ generated by all elements of $G^+$ such that $[x]$ belongs
to $S_{\omega}^{(\omega)}$. However, it follows from
Corollary~\ref{C:CXAl1} that this is not possible. So we have obtained the
following negative result: \emph{There exists a distributive \jzus\
without descending $\omega_1$-chains that cannot be isomorphic to
$\nabla(G^+)$ for an interpolation group $G$}.

This implies, in turn, the following negative result: \emph{$\URPSR$
is not sufficient to characterize all distributive \jzus s of the form
$\nabla(G^+)$ for $G$ an interpolation group}.
\end{remark}

\section{Open problems}\label{S:Pbs}

Our most intriguing open problem is related to the following known results:
\begin{itemize}\em
\item[---] Every \jzs\ of the form $\varinjlim_{n<\omega}D_n$, with
all $D_n$-s being distributive lattices with zero and all transition maps
being \jzh s, is isomorphic to the semilattice
of all compact congruences of some relatively complemented lattice with zero
\pup{see \cite{Wehr}}.

\item[---] Every distributive lattice with zero is isomorphic to
$\nabla(G^+)$, for some dimension group $G$ \pup{see \cite{GoWe1}}.
\end{itemize}

Is it possible to unify these results? We can, for example, ask the
following.

\begin{problem}\label{Pb:UnDn}
Let $S=\varinjlim_{n<\omega}D_n$, with all $D_n$-s being distributive
lattices with zero and all transition maps being \jzh s.
Does there exist a dimension group $G$ such that $S\cong\nabla(G^+)$?
\end{problem}

The uniform refinement property $\URPSR$ is of no help to solve
Problem~\ref{Pb:UnDn} negatively: indeed, by Corollary~\ref{C:UnDn}, $S$
does satisfy $\URPSR$.

\begin{problem}\label{Pb:IdSSIG}
Let $M$ be a strongly separative refinement monoid. Does there exist an
interpolation group $G$ such that $\nabla(M)\cong\nabla(G^+)$?
\end{problem}

As in \cite{Ruzi1}, we say that a \poag\ $G$ is \emph{weakly Archimedean},
if for all $a$, $b\in G^+$, if $na\leq b$ holds for all $n\in\NN$, then
$a=0$. By using one of the main results in \cite{TuWe1}, P.
R\r{u}\v{z}i\v{c}ka proves in \cite{Ruzi1} the following result: \emph{Every
countable distributive \jzs\ is isomorphic to $\nabla(G^+)$, for some weakly
Archimedean dimension group $G$}, see \cite[Theorem~3.1]{Ruzi1}.

We say that a \poag\ $G$ is \emph{Archimedean}, if for all $a$,
$b\in G$, if $na\leq b$ for all $n\in\NN$, then $a\leq 0$.

\begin{problem}\label{Pb:ArchLift}
For a countable distributive \jzs\ $S$, does there exist an Archimedean
dimension group $G$ such that $\nabla(G^+)\cong S$?
\end{problem}

Another problem, inspired by Bergman's Theorem,
Corollary~\ref{C:LocFinLatt}, and some open questions in \cite{TuWe2},
is the following.

\begin{problem}\label{Pb:LocFin}
For a countable distributive semilattice $S$, does there exist a modular
lattice $K$, generating a locally finite variety, such that
$\Conc K\cong S$?
\end{problem}

The statement that $K$ generates a locally finite variety is stronger than
the mere local finiteness of $K$. It is equivalent to saying that for every
$n\in\NN$, the cardinalities of all $n$-generated sublattices of $K$ are
bounded by a positive integer.

Our next problem asks about lifting not semilattices, but \emph{diagrams} of
semilattices. We say that a diagram of semilattices, indexed by a partially
ordered set~$I$, is \emph{finite} (resp., \emph{countable},
\emph{dismantlable}), if $I$ is finite (resp., countable, dismantlable).
It is proved in \cite{TuWe1} that every finite dismantlable diagram of
finite Boolean semilattices can be lifted, with respect to the $\nabla$
functor, by a diagram of (positive cones of) dimension groups.

\begin{problem}\label{Pb:DismLift}
Can every countable dismantlable diagram of countable distributive \jzs s
be lifted, with respect to the $\nabla$ functor, by a diagram of (positive
cones of) dimension groups?
\end{problem}

\section*{Acknowledgment}
I thank the anonymous referee for his thoughtful report, which
brought additional life and openness to the topics discussed in the paper.


\begin{thebibliography}{99}
\bibitem{AdGo}
K.\,V. Adaricheva and V.\,A. Gorbunov,
\emph{On lower bounded lattices},
Algebra Universalis \textbf{46} (2001), 203--213.

\bibitem{AraS}
P. Ara,
\emph{Stability properties of exchange rings},
International Symposium on Ring Theory (Kyongju, 1999),
Trends Math., Birkh\"auser Verlag, Boston, MA, 2001, 23--42.

\bibitem{AGPO}
P. Ara, K.\,R. Goodearl, K.\,C. O'Meara, and E. Pardo,
\emph{Separative cancellation for projective modules over exchange rings},
Israel J. Math. \textbf{105} (1998), 105--137.

\bibitem{Berg86}
G.\,M. Bergman,
\emph{Von Neumann regular rings with tailor-made ideal lattices},
Unpublished note (26 October 1986).

\bibitem{Day79}
A. Day,
\emph{Characterization of finite lattices that are bounded homomorphic
images of sublattices of free lattices},
Canad. J. Math. \textbf{31} (1979), 69--78.

\bibitem{EHS80}
E.\,G. Effros, D.\,E. Handelman, and C.-L. Shen,
\emph{Dimension groups and their affine representations},
Amer. J. Math. \textbf{102} (1980), no.~2,
385--407.

\bibitem{FJN}
R. Freese, J. Je\v{z}ek, and J.\,B. Nation,
``Free Lattices'',
Mathematical Surveys and Monographs \textbf{42}, Amer. Math. Soc.,
Providence, 1995. viii+293~p.

\bibitem{Gvnrr}
K.\,R. Goodearl,
``Von Neumann Regular Rings'', Pitman, London 1979;
Second Ed. Krieger, Malabar, Fl., 1991.

\bibitem{Gpoag}
\bysame,
``Partially Ordered Abelian Groups with Interpolation'',
Math. Surveys and Monographs \textbf{20}, Amer. Math. Soc.,
Providence, 1986.

\bibitem{GoHa86}
K.\,R. Goodearl and D.\,E. Handelman,
\emph{Tensor products of dimension groups and $K_0$ of
unit-regular rings}, Canad. J. Math. \textbf{38}, no.~3 (1986), 633--658.

\bibitem{GoWe1}
K.\,R. Goodearl and F. Wehrung,
\emph{Representations of distributive semilattices in ideal lattices of
various algebraic structures}, Algebra Universalis~\textbf{45} (2001),
71--102.

\bibitem{Grat98}
G. Gr\"atzer,
``General Lattice Theory. Second edition'', new appendices by the
author with B.\,A. Davey, R. Freese, B. Ganter,
M. Greferath, P. Jipsen, H.\,A. Priestley, H. Rose, E.\,T. Schmidt,
S.\,E. Schmidt, F. Wehrung, and R. Wille. Birkh\"auser Verlag, Basel,
1998. xx+663~p.

\bibitem{GLWe}
G.~Gr{\"a}tzer, H.~Lakser, and F.~Wehrung,
\emph{Congruence amalgamation of lattices},
Acta Sci. Math. (Szeged) \textbf{66} (2000), 339--358.

\bibitem{GrWe}
G. Gr\"atzer and F. Wehrung,
\emph{On the number of join-irreducibles in a
congruence representation of a finite distributive lattice},
Algebra Universalis \textbf{49} (2003), 165--178.

\bibitem{JoKi}
B. J\'onsson and J.\,E. Kiefer,
\emph{Finite sublattices of a free lattice},
Canad. J. Math. \textbf{14} (1962), 487--497.

\bibitem{CMor2}
C. Moreira dos Santos,
\emph{A refinement monoid whose maximal antisymmetric quotient is not a
refinement monoid}, Semigroup Forum \textbf{65}, no.~2 (2002), 249--263.

\bibitem{Pard95}
E. Pardo,
\emph{Monoides de refinament i anells d'intercanvi}, PhD. Thesis,
Universitat Aut\`onoma de Barcelona, 1995.

\bibitem{Ruzi1}
P. R\r{u}\v{z}i\v{c}ka,
\emph{A distributive semilattice not isomorphic to the maximal semilattice
quotient of the positive cone of any dimension group}, J. Algebra
\textbf{268}, no.~1 (2003), 290--300.

\bibitem{TuWe1}
J. T\r{u}ma and F. Wehrung,
\emph{Liftings of diagrams of semilattices by diagrams of dimension
groups}, Proc. London Math. Soc. \textbf{87}, no.~3 (2003), 1--28.

\bibitem{TuWe2}
\bysame,
\emph{A survey of recent results on congruence lattices of lattices},
Algebra Universalis \textbf{48}, no.~4 (2002), 439--471.

\bibitem{Wehr98a}
F. Wehrung,
\emph{Non-measurability properties of interpolation vector spaces},
Israel J. Math. \textbf{103} (1998), 177--206.

\bibitem{Wehr98b}
\bysame,
\emph{The dimension monoid of a lattice},
Algebra Universalis \textbf{40}, no.~3 (1998), 247--411.

\bibitem{Wehr99}
\bysame,
\emph{A uniform refinement property for congruence lat\-tices},
Proc. Amer. Math. Soc. \textbf{127}, no.~2 (1999), 363--370.

\bibitem{Wehr00}
\bysame,
\emph{Representation of algebraic distributive lattices with $\aleph_1$ compact
elements as ideal lattices of regular rings},
Publ. Mat. (Barcelona) \textbf{44}, no.~2 (2000), 419--435.

\bibitem{Wehr02}
\bysame,
\emph{From join-irreducibles to dimension theory for lattices with chain
conditions}, J.~Algebra Appl. \textbf{1}, no.~2
(2002), 215--242.

\bibitem{Wehr}
\bysame,
\emph{Forcing extensions of partial lattices},
J. Algebra \textbf{262}, no.~1 (2003), 127--193.

\end{thebibliography}
\end{document}